\magnification=\magstep1
\overfullrule=0pt
\tolerance=400
\font\rms=cmr9
\font\bfl=cmbx10 scaled \magstep1
\font\mb=msbm10
\font\mbs=msbm7
\font\mbss=msbm5
\newfam\mbfam
\textfont\mbfam=\mb
\scriptfont\mbfam=\mbs
\scriptscriptfont\mbfam=\mbss
\def\mb{\fam\mbfam}
\def\npar{\relax\par\noindent}
\def\nasp{\spacefactor=1000}
\def\\{\ifhmode\hfil\break\else\ifvmode\vskip\baselineskip\fi\fi}
\def\vbreak#1#2{\par\ifdim\lastskip<#2\removelastskip\penalty#1\vskip#2\fi}
\def\vbreakn#1#2{\vbreak{#1}{#2}\npar}
\def\title#1{\vbreak{-9000}{1cm}\centerline{\bfl#1}}
\def\author#1{\vbreak{9000}{5mm}\centerline{#1}}
\def\titledatefix#1{\vbreak{9000}{5mm}\centerline{\it #1}}
\long\def\abstract#1{\vbreak{9000}{5mm}\centerline{\bf Abstract}
 \kern3mm{#1}\vbreak{-9000}{4mm}}
\def\support#1{\footnote{}{{\rms#1\hfil}}}
\def\acknowledgment#1{\vbreakn{-5000}{2mm}{\bf Acknowledgment.} #1}
\def\section#1#2{\vbreak{-9500}{6mm}\centerline{\bf #1. #2}\vbreak{9900}{3mm}}

\def\equ#1{$$\vcenter{\equalign{#1}}$$}
\def\equalign#1{\let\\=\cr\let\-=\hfill
\ialign{&\hfil$\dsp ##$\hfil\cr#1\crcr}}
\def\lequ#1{$$\vcenter{\lequalign{#1}}$$}
\def\lequalign#1{\let\\=\cr\let\-=\hfill
\ialign{&\hbox to \hsize{$\dsp ##$\hfil}\cr#1\crcr}}
\def\state#1#2#3{\vbreakn{-5000}{2mm}{\bf #1 #2.\nasp\ }#3\begingroup}
\def\nstate#1{\vbreakn{-5000}{2mm}{\bf #1.\nasp\ }\begingroup}

\def\endstate{\par\endgroup\vbreak{-8000}{2mm}}
\def\theorem#1#2#3{\state{Theorem}{#2}{#3}\it}
\def\endtheorem{\endstate}
\def\lemma#1#2#3{\state{Lemma}{#2}{#3}\it}
\def\endlemma{\endstate}
\def\proposition#1#2#3{\state{Proposition}{#2}{#3}\it}
\def\endproposition{\endstate}
\def\corollary#1#2#3{\state{Corollary}{#2}{#3}\it}
\def\endcorollary{\endstate}
\def\proof{\nstate{Proof}}
\def\endproof{\endprr\endstate}
\def\enprule{\vrule height1mm depth1mm width2mm}
\def\endprr{\discretionary{}{\kern\hsize}{\kern 3ex}\hbox to 0pt{\hss\enprule}}
\def\rfr#1#2{#1~#2}
\def\bibliography{\vbreak{-9000}{6mm}
\line{\bf Bibliography\hfil}\kern2mm
\begingroup\parskip=0pt\parindent=0pt\frenchspacing\rms}
\def\endbibliography{\par\endgroup}
\def\bibitem#1#2{\vbreak{-9500}{1mm}\hangindent=11mm \hangafter=1
\noindent\hbox to 0pt{[#1]\hss}\hskip\hangindent #2}
\def\bibart#1 a:#2 t:#3 j:#4 n:#5 y:#6 p:#7 *{%
\bibitem{#1}{#2, #3, {\it #4\/} {\bf #5} (#6), #7.}}
\def\bibook#1 a:#2 t:#3 i:#4 *{%
\bibitem{#1}{#2, {\it #3\/}, #4.}}
\def\bibartp#1 a:#2 t:#3 *{\bibitem{#1}{#2, #3, in preparation.}}
\def\bibartpr#1 a:#2 t:#3 *{\bibitem{#1}{#2, #3, preprint.}}
\def\bibartn#1 a:#2 t:#3 j:#4 *{\bibitem{#1}{#2, #3, #4.}}
\def\no{no. }
\def\brfr#1{[#1]}
\def\sdup#1#2{{\scr #1 \atop \scr #2}}

\def\rest#1{\raise-2pt\hbox{$|_{#1}$}}
\def\frac#1#2{{#1\over #2}}

\def\comp{\mathord{\hbox{$\scr\circ$}}}

\def\dsc{\discretionary{}{}{}}
\def\R{{\mb R}}
\def\C{{\mb C}}
\def\Z{{\mb Z}}
\def\N{{\mb N}}

\def\mod{\mathop{\rm mod}}
\def\ed#1{1_{#1}}
\let\dsp=\displaystyle

\let\ovr=\overline
\let\scr=\scriptstyle

\let\sln=\subset
\let\sle=\subseteq

\let\sm=\setminus
\let\col=\colon

\let\ld=\ldots

\let\cd=\cdot

\let\ra=\longrightarrow
\let\ras=\rightarrow
\let\tens=\otimes
\let\alf=\alpha

\let\Gam=\Gamma
\let\gam=\gamma
\let\del=\delta

\let\phi=\varphi
\let\sig=\sigma
\let\eps=\varepsilon
\let\lam=\lambda
\let\om=\omega
\let\ro=\rho
\def\hX{\widehat{X}}
\def\hG{\widehat{G}}

\def\hx{\hat{x}}
\def\hf{\hat{f}}
\def\tf{\tilde{f}}
\def\tY{\widetilde{Y}}
\def\B{{\cal B}}
\def\F{{\cal F}}
\def\cZ{{\cal Z}}
\def\cN{{\cal N}}
\def\cNb{\cN^{o}}
\def\cM{{\cal M}}
\def\cMb{\cM^{o}}
\def\cP{{\cal P}}
\def\cPb{\cP^{o}}

\def\cD{{\cal D}}
\def\cV{{\cal V}}
\def\cL{{\cal L}}
\def\cQ{{\cal Q}}
\def\cX{{\cal X}}
\def\Gc{G^{c}}
\def\Hc{H^{c}}

\def\tX{\widetilde{X}}
\def\tG{\widetilde{G}}
\def\ta{\tilde{a}}
\def\tb{\tilde{b}}
\def\tc{\tilde{c}}
\def\tx{\tilde{x}}
\def\ty{\tilde{y}}
\def\tu{\tilde{u}}
\def\ha{\hat{a}}
\def\hh{\hat{h}}
\def\hpi{\hat{\pi}}
\def\tpi{\tilde{\pi}}
\let\Lam=\Lambda
\let\Om=\Omega
\let\ophi=\oldphi
\let\ups=\upsilon
\let\ro=\ro
\def\tro{\tilde{\ro}}
\def\T{{\mb T}}
\let\dX=\partial
\def\dist{\mathop{\hbox{\rm dist}}}
\def\lfa#1#2{\mathord{\hbox{\vtop{\kern-6pt\hbox{$#2$}}$\setminus$$#1$}}}
\def\Li{l^{\infty}}
\def\supp{\mathop{\hbox{\rm supp}}}
\title{Nilsequences, null-sequences, and multiple correlation sequences}
\author{A. Leibman}
\support{Supported by NSF grants DMS-0901106 and DMS-1162073.}
\kern3mm
\centerline{\vbox{\hsize=4.5cm\rms\noindent
Department of Mathematics\\
The Ohio State University\\
Columbus, OH 43210, USA\\
e-mail: leibman@math.ohio-state.edu}}
\titledatefix{April 3, 2013}
\abstract{
A ($d$-parameter) {\it basic nilsequence\/} 
is a sequence of the form $\psi(n)=f(a^{n}x)$, $n\in\Z^{d}$,
where $x$ is a point of a compact nilmanifold $X$,
$a$ is a translation on $X$,
and $f\in C(X)$;
a {\it nilsequence\/} is a uniform limit of basic nilsequences.
If $X=G/\Gam$ be a compact nilmanifold,
$Y$ is a subnilmanifold of $X$,
$(g(n))_{n\in\Z^{d}}$ is a polynomial sequence in $G$, and $f\in C(X)$,
we show that the sequence $\phi(n)=\int_{g(n)Y}f$
is the sum of a basic nilsequence 
and a sequence that converges to zero in uniform density (a null-sequence).
We also show that an integral of a family of nilsequences is a nilsequence plus a null-sequence.
We deduce that
for any invertible finite measure preserving system $(W,\B,\mu,T)$,
polynomials $p_{1},\ld,p_{k}\col\Z^{d}\ra\Z$,
and sets $A_{1},\ld,A_{k}\in\B$,
the sequence 
$\phi(n)=\mu(T^{p_{1}(n)}A_{1}\cap\ld\cap T^{p_{k}(n)}A_{k})$, $n\in\Z^{d}$,
is the sum of a nilsequence and a null-sequence.}
\kern3mm
\section{0}{Introduction}
Throughout the whole paper we will deal with ``multiparameter sequences'', --
we fix $d\in\N$
and under ``a sequence'' will usually understand ``a two-sided $d$-parameter sequence'',
that is, a mapping with domain $\Z^{d}$.

A (compact) $r$-step {\it nilmanifold\/} $X$ is a factor space $G/\Gam$,
where $G$ is an $r$-step nilpotent (not necessarily connected) Lie group
and $\Gam$ is a discrete co-compact subgroup of $G$.
Elements of $G$ act on $X$ by {\it translations};
an ($r$-step) {\it nilsystem\/} is an ($r$-step) nilmanifold $X=G/\Gam$
with a translation $a\in G$ on it.

{\it A basic $r$-step nilsequence\/} 
is a sequence of the form $\psi(n)=f(\eta(n)x)$, $n\in\Z^{d}$,
where $x$ is a point of an $r$-step nilmanifold $X=G/\Gam$,
$\eta$ is a homomorphism $\Z^{d}\ra G$, and $f\in C(X)$;
{\it an $r$-step nilsequence\/} is a uniform limit of basic $r$-step nilsequences.
The algebra of nilsequences is a natural generalization
of Weyl's algebra of almost periodic sequences,
which are just 1-step nilsequences.
An ``inner'' characterization of nilsequences, in terms of their properties,
is obtained in \brfr{HKM};
see also \brfr{HK2}.

The term ``nilsequence'' was introduced in \brfr{BHK},
where it was proved that for any ergodic finite measure preserving system $(W,\B,\mu,T)$,
positive integer $k$, and sets $A_{1},\ld,A_{k}\in\B$
the {\it multiple correlation sequence\/}
$\phi(n)=\mu\bigl(T^{n}A_{1}\cap\ld\cap T^{kn}A_{k}\bigr)$, $n\in\N$,
is a sum of a nilsequence 
and of a sequence that tends to zero in uniform density in $\Z^{d}$,
{\it a null-sequence}.
Our goal in this paper is to generalize this result
to multiparameter polynomial multiple correlation sequences and general (non-ergodic) systems.
We prove (in \rfr{Section}{5}):
\theorem{P-ICor}{0.1}{}
Let $(W,\B,\mu,T)$ be an invertible measure preserving system with $\mu(W)<\infty$,
let $p_{1},\ld,p_{k}$ be polynomials $\Z^{d}\ra\Z$,
and let $A_{1},\ld,A_{k}\in\B$.
Then the ``multiple polynomial correlation sequence'' 
$\phi(n)=\mu\bigl(T^{p_{1}(n)}A_{1}\cap\ld\cap T^{p_{k}(n)}A_{k}\bigr)$, $n\in\Z^{d}$,
is a sum of a nilsequence and a null-sequence.
\endtheorem
\npar
(In \brfr{L6} this theorem was proved in the case $d=1$ and ergodic $T$.)

Based on the theory of nil-factors developed in \brfr{HK1} and, independently, in \brfr{Z},
it is shown in \brfr{L3} that nilsystems are {\it characteristic\/}
for multiple polynomial correlation sequences induced by ergodic systems,
in the sense that, up to a null-sequence and an arbitrarily small sequence,
any such correlation sequence comes from a nilsystem.
This reduces the problem of studying ``ergodic'' multiple polynomial correlation sequences
to nilsystems.

Let $X=G/\Gam$ be a connected nilmanifold,
let $Y$ be a connected subnilmanifold of $X$,
and let $g$ be a polynomial sequence in $G$,
that is, a mapping $\Z^{d}\ra G$
of the form $g(n)=a_{1}^{p_{1}(n)}\ld a_{r}^{p_{r}(n)}$, $n\in\Z^{d}$,
where $a_{1},\ld,a_{r}\in G$
and $p_{1},\ld,p_{r}$ are polynomials $\Z^{d}\ra\Z$.
We investigate (in \rfr{Section}{3}) the behavior of the sequence $g(n)Y$ of subnilmanifolds of $X$:
we show that there is a subnilmanifold $Z$ of $X$, containing $Y$,
such that the sequence $g(n)$ only shifts $Z$ along $X$, without distorting it,
whereas, outside of a null-set of $n\in\Z^{d}$,
$g(n)Y$ becomes more and more ``dense'' in $g(n)Z$:
\proposition{P-IYtoZ}{0.2}{}
Assume {\rm(as we can)} that the orbit $g(n)Y$, $n\in\Z^{d}$, is dense in $X$,
and let $Z$ be the normal closure {\rm(in the algebraic sense; see below)} of $Y$ in $X$.
Then for any $f\in C(X)$, 
the sequence $\lam(n)=\int_{g(n)Y}f-\int_{g(n)Z}f$, $n\in\Z^{d}$, is a null-sequence.
\endproposition
We have $\int_{g(n)Z}f=g(n)\hf(g(n)e)$, $n\in\Z^{d}$, where $\hf=E(f|X/Z)$ and $e=Z/Z\in X/Z$.
(Here and below, 
$E(f|X')$ stands for the conditional expectation of a function $f\in L^{1}(X)$
with respect to a factor $X'$ of $X$.)
So, the sequence $\int_{g(n)Z}f$ is a basic nilseqeunce, and we obtain:
\theorem{P-IYnils}{0.3}{}
For any $f\in C(X)$
the sequence $\phi(n)=\int_{g(n)Y}f$, $n\in\Z^{d}$,
is the sum of a basic nilsequence and a null-sequence.
\endtheorem
Applying this result to the diagonal $Y$ of the power $X^{k}$ of the nilmanifold $X$,
the polynomial sequence $g(n)=(a^{p_{1}(n)},\ld,a^{p_{k}(n)})$, $n\in\Z^{d}$, in $G^{k}$,
and the function $f=1_{A_{1}}\otimes\ld\otimes1_{A_{k}}$,
we obtain \rfr{Theorem}{0.1} in the ergodic case.

Our next step (\rfr{Section}{4}) is to extend this result to the case of a non-ergodic $T$.
Using the ergodic decomposition $W\ra\Om$ of $T$
we obtain a measurable mapping from $\Om$ to the space of nilsequences--plus--null-sequences,
which we then have to integrate over $\Om$.
The integral of a family of null-sequences is a null-sequence, and creates no trouble.
As for nilsequences, when we integrate them we arrive at the following problem:
if $X=G/\Gam$ is a nilmanifold, with $\pi\col G\ra X$ being the factor mapping, 
and $\ro(a)$, $a\in G$, is a finite Borel measure on $G$,
what is the limiting behavior of the measures $\pi_{*}(\ro(a^{n}))$ on $X$?
(This is the question corresponding to the case $d=1$;
for $d\geq 2$ it is slightly more complicated.)
We show that this sequence of measures tends
to a linear combination of Haar measures on (countably many) subnilmanifolds of $X$,
which are normal (and so travel, without distortion) in the closure of their orbits,
and we again obtain:
\proposition{P-IlimMes}{0.4}{}
For any $f\in C(X)$, the sequence $\phi(n)=\int_{G}f(\pi(a^{n}))\,d\ro(a)$, $n\in\Z$,
is a sum of a basic nilsequence and a null-sequence.
\endproposition
\npar
(This proposition is a ``nilpotent'' extension of the following classical fact:
if $\ro$ is a finite Borel measure on the 1-dimensional torus $\T$,
then its Fourier transform
$\phi(n)=\int_{\T}e^{-2\pi inx}d\ro(x)$ is the sum
of an almost periodic sequence
(a 1-step nilsequence; it corresponds to the atomic part of $\ro$)
and a null-sequence
(that corresponds to the non-atomic part of $\ro$).)

As a corollary, we obtain the remaining ingredient of the proof of \rfr{Theorem}{0.1}:
\theorem{P-IInt}{0.5}{}
Let $\Om$ be a measure space
and let $\phi_{\om}$, $\om\in\Om$, be an integrable family of nilsequences;
then the sequence $\phi(n)=\int_{\Om}\phi_{\om}(n)$
is a sum of a nilsequence and a null-sequence.
\endtheorem

Let us also mention {\it generalized {\rm(or {\it bracket})} polynomials}, --
the functions constructed from ordinary polynomials
using the operations of addition, multiplication, and taking the integer part, $[\cd]$.
(For example, $p_{1}[p_{2}[p_{3}]+p_{4}]$, where $p_{i}$ are ordinary polynomials,
is a generalized polynomial.)
Generalized polynomials (gps) appear quite often
(for example, the fractional part, and the distance to the nearest integer, 
of an ordinary polynomial are gps);
they were systematically studied in \brfr{H{\aa}1}, \brfr{H{\aa}2}, \brfr{BL}, and \brfr{L7}.
Because of their simple definition, gps are nice objects to deal with.
On the other hand, similarly to nilsequences, gps come from nilsystems:
bounded gps (on $\Z^{d}$, in our case)
are exactly the sequences of the form $h(g(n)x)$, $n\in\Z^{d}$,
where $h$ is a piecewise polynomial function on a nilmanifold $X=G/\Gam$, $x\in X$,
and $g$ is a polynomial sequence in $G$
(see \brfr{BL} or \brfr{L7}).
Since any continuous function is uniformly approximable by piecewise polynomial functions
(this follows by an application of the Weierstrass theorem in the fundamental domain of $X$),
nilsequences are uniformly approximable by generalized polynomials.
We obtain as a corollary that any multiple polynomial correlation sequence
is, up to a null-sequence, uniformly approximable by generalized polynomials.
\acknowledgment
I am grateful to Vitaly Bergelson for his interest and advice.
\section{1}{Nilmanifolds}

In this section we collect the facts about nilmanifolds that we will need below;
details and proofs can be found in \brfr{M}, \brfr{L1}, \brfr{L2}, \brfr{L4}.

Throughout the paper,
$X=G/\Gam$ will be a compact nilmanifold,
where $G$ is a nilpotent Lie group 
and $\Gam$ is a discrete subgroup of $G$,
and $\pi$ will denote the factor mapping $G\ra X$.
By $\ed{X}$ we will denote the point $\pi(\ed{G})$ of $X$.
By $\mu_{X}$ we will denote the normalized Haar measure on $X$.

By $\Gc$ we will denote the identity component of $G$.
Note that if $G$ is disconnected,
$X$ can be interpreted as a nilmanifold, $X=G'/\Gam'$, in different ways:
if, for example, $X$ is connected, then also $X=\Gc/(\Gam\cap\Gc)$.
If $X$ is connected
and we study the action on $X$ of a sequence $g(n)$, $g\col\Z^{d}\ra G$,
we may always assume that $G$ is generated by $\Gc$
and the range $g(\Z^{d})$ of $g$.
Thus, we may (and will) assume that the group $G/\Gc$ is finitely generated.

Every nilpotent Lie group $G$ is a factor of a torsion free nilpotent Lie group.
(As such, a suitable ``free nilpotent Lie group'' $F$ can be taken.
If $\Gc$ has $k_{1}$ generators, $G/\Gc$ has $k_{2}$ generators,
and $G$ is $r$-step nilpotent, then $F=\F/\F_{r+1}$,
where $\F$ is the free product
of $k_{1}$ copies of $\R$ and $k_{2}$ copies of $\Z$,
and $\F_{r+1}$ is the $(r+1)$st term of the lower central series of $\F$.)
Thus, we may always assume that $G$ is torsion-free.
The identity component $\Gc$ of $G$ is then an exponential Lie group,
and for every element $a\in\Gc$ 
there exists a (unique) one-parametric subgroup $a^{t}$ such that $a^{1}=a$.

If $G$ is torsion free,
it possesses a {\it Malcev basis\/} compatible with $\Gam$,
which is a finite set $\{e_{1},\ld,e_{k}\}$ of elements of $\Gam$,
with $e_{1},\ld,e_{k_{1}}\in\Gc$ and $e_{k_{1}+1},\ld,e_{k}\not\in\Gc$,
such that every element $a\in G$ can be uniquely written in the form
$a=e_{1}^{u_{1}}\ld e_{k}^{u_{k}}$
with $u_{1},\ld,u_{k_{1}}\in\R$ and $u_{k_{1}+1},\ld,u_{k}\in\Z$,
and with $a\in\Gam$ iff $u_{1},\ld,u_{k}\in\Z$;
we call $u_{1},\ld,u_{k}$ {\it the coordinates\/} of $a$.
Thus, Malcev coordinates define a homeomorphism 
$G\simeq\R^{k_{1}}\times\Z^{k-k_{1}}$, $a\leftrightarrow(u_{1},\ld,u_{k})$,
which maps $\Gam$ onto $\Z^{k}$.

The multiplication in $G$ is defined by the (finite) multiplication table 
for the Malcev basis of $G$, whose entries are integers;
it follows that there are only countably many non-isomorphic nilpotent Lie groups
with cocompact discrete subgroups,
and countably many non-isomorphic compact nilmanifolds.

Let $X$ be connected. 
Then, under the identification $\Gc\leftrightarrow\R^{k_{1}}$,
the cube $[0,1)^{k_{1}}$ is the fundamental domain of $X$.
We will call the closed cube $Q=[0,1]^{k_{1}}$ 
{\it the fundamental cube of $X$ in $\Gc$}
and simetimes identify $X$ with $Q$.
When $X$ is identified with its fundamental cube $Q$,
the measure $\mu_{X}$ corresponds to the standard Lebesgue measure $\mu_{Q}$ on $Q$.

In Malcev coordinates,
multiplication in $G$ is a polynomial operation:
there are polynomials $q_{1},\ld,q_{k}$ 
in $2k$ variables with rational coefficients
such that for $a=e_{1}^{u_{1}}\ld e_{k}^{u_{k}}$ 
and $b=e_{1}^{v_{1}}\ld e_{k}^{v_{k}}$
we have $ab=e_{1}^{q_{1}(u_{1},v_{1},\ld,u_{k},v_{k})}\ld
e_{k}^{q_{k}(u_{1},v_{1},\ld,u_{k},v_{k})}$.
This implies that ``life is polynomial'' in nilpotent Lie groups:
in coordinates, homomorphisms between these groups are polynomial mappings,
and connected closed subgroups of such groups 
are images of polynomial mappings
and are defined by systems of polynomial equations.

{\it A subnilmanifold\/} $Y$ of $X$ is a closed subset of the form $Y=Hx$,
where $H$ is a closed subgroup of $G$ and $x\in X$.
For a closed subgroup $H$ of $G$,
the set $\pi(H)=H\ed{X}$ is closed (and so, is a subnilmanifold)
iff the subgroup $\Gam\cap H$ is co-compact in $H$;
we will call the subgroup $H$ with this property {\it rational}.
Any subnilmanifold $Y$ of $X$ has the form $\pi(aH)=a\pi(H)$,
where $H$ is a closed rational subgroup of $G$.

If $Y$ is a subnilmanifold of $X$ with $\ed{X}\in Y$,
then $Y=\pi(H)=H\ed{X}$ for some closed subgroup $H$ of $G$.
$H$ may not be the minimal subgroup with this property:
if $Y$ is connected,
then the identity component $\Hc$ of $H$ also satisfies $\pi(\Hc)=Y$.

The intersection of two subnilmanifolds is a finite disjoint union of subnilmanifolds.

Given a subnilmanifold $Y$ of $X$,
by $\mu_{Y}$ we will denote the normalized Haar measure on $Y$.
Translations of subnilmanifods are measure preserving:
we have $a_{*}\mu_{Y}=\mu_{aY}$ for all $a\in G$.

Let $Z$ be a subnilmanifold of $X$, $Z=Lx$,
where $L$ is a closed subgroup of $G$.
We say that $Z$ is {\it normal\/} if $L$ is normal.
In this case the nilmanifold $\hX=X/Z=G/(L\Gam)$ is defined,
and $X$ splits into a disjoint union of fibers
of the factor mapping $X\ra\hX$.
(Note that if $L$ is normal in $\Gc$ only,
then the factor $X/Z=\Gc/(L\Gam)$ is also defined,
but the elements of $G\sm\Gc$ do not act on it.)

One can show that a subgroup $L$ is normal
iff $\gam L\gam^{-1}=L$ for all $\gam\in\Gam$;
hence, $Z=\pi(L)$ is normal iff $\gam Z=Z$ for all $\gam\in\Gam$.

If $H$ is a closed rational subgroup of $G$
then its normal closure $L$ 
(the minimal normal subgroup of $G$ containing $H$)
is also closed and rational,
thus $Z=\pi(L)$ is a subnilmanifold of $X$.
We will call $Z$ {\it the normal closure\/} 
of the subnilmanifold $Y=\pi(H)$.
If $L$ is normal then the identity component of $L$ is also normal;
this implies that the normal closure of a connected subnilmanifold
is connected.

If $X$ is connected,
{\it the maximal factor-torus\/} of $X$ is the torus $\lfa{X}{[\Gc,\Gc]}$,
and {\it the nil-maximal factor-torus\/} is $\lfa{X}{[G,G]}$.
The nil-maximal factor-torus is a factor of the maximal one.

If $\eta\col\Z^{d}\ra G$ is a homomorphism,
then for any point $x\in X$ the closure of the orbit $\eta(\Z^{d})x$ of $x$ in $X$
is a subnilmanifold $V$ of $X$ (not necessarily connected),
and the sequence $\eta(n)x$, $n\in\Z^{d}$, is well distributed in $V$.
(This means that for any function $f\in C(V)$ and any F{\o}lner sequence $(\Phi_{N})$ in $\Z^{d}$,
$\lim_{N\ras\infty}\frac{1}{|\Phi_{N}|}\sum_{n\in\Phi_{N}}f(\eta(n)x)=\int_{Y}f\,d\mu_{V}$.)
If $X$ is connected,
the sequence $\eta(n)x$, $n\in\Z^{d}$, is dense, and so, well distributed in $X$
iff the image of this sequence is dense in the nil-maximal factor-torus of $X$.
All the same is true for the orbit of any subnilmanifold $Y$ of $X$:
the closure of $\bigcup_{n\in\Z^{d}}\eta(n)Y$ is a subnilmanifold $W$ of $X$;
the sequence $\eta(n)Y$, $n\in\Z^{d}$, is well distributed in $W$
(this means that for any function $f\in C(W)$ and any F{\o}lner sequence $(\Phi_{N})$ in $\Z^{d}$,
$\lim_{N\ras\infty}\frac{1}{|\Phi_{N}|}\sum_{n\in\Phi_{N}}\int_{\eta(n)Y}f(x)\,d\mu_{\eta(n)Y}
=\int_{Y}f\,d\mu_{V}$);
and, in the case $X$ is connected, 
the sequence $\eta(n)Y$ is well distributed in $X$
iff its image is dense in the nil-maximal factor-torus of $X$.

{\it A polynomial sequence\/} in $G$ is a sequence of the form
$g(n)=a_{1}^{p_{1}(n)}\kern-3mm\ld a_{k}^{p_{k}(n)}$, $n\in\Z^{d}$,
where $a_{1},\ld,a_{k}\in G$ and $p_{1},\ld,p_{k}$ are polynomials $\Z^{d}\ra\Z$.
Let $g$ be a polynomial sequence in $G$ and let $x\in X$.
Then the closure $V$ of the orbit $g(\Z^{d})x$ 
is a finite disjoint union of connected subnilmanifolds of $X$,
and $g(n)x$ visits these subnilmanifolds periodically:
there exists $l\in\N$ such that for any $i\in\Z^{d}$, 
all the elements $g(lm+i)x$, $m\in\Z^{d}$, belong to the same connected component of $V$.
If $V$ is connected, then the sequence $g(n)x$, $n\in\Z^{d}$, is well distributed in $V$.
In the case $X$ is connected, 
the sequence $g(n)x$, $n\in\Z^{d}$, is dense, and so, well distributed in $X$
iff the image of this sequence is dense in the maximal factor-torus of $X$.
All the same is true for the orbit $g(Z^{d})Y$ 
of any connected subnilmanifold $Y$ of $X$ under the action of $g$:
its closure $W$ is a finite disjoint union of connected subnilmanifolds of $X$,
visited periodically;
if $W$ is connected, then the sequence $g(n)Y$, $n\in\Z^{d}$, is well distributed in $W$;
and, if $X$ is connected, the sequence $g(n)Y$ is well distributed in $X$
iff its image is dense in the maximal factor-torus of $X$.

The following proposition, which is a corollary 
(of a special case) of the result obtained in \brfr{GT},
says that ``almost every'' subnilmanifold of $X$ 
is distributed in $X$ ``quite uniformly''.
(See Appendix in \brfr{L6} for details.)
\proposition{P-myTao}{1.1}{}
Let $X$ be connected.
For any $C>0$ and any $\eps>0$ 
there are finitely many subnilmanifolds $V_{1},\ld,V_{r}$ of $X$,
connected and containing $\ed{X}$,
such that for any connected subnilmanifold $Y$ of $X$ with $\ed{X}\in Y$, 
if $Y\not\sle V_{i}$ for all $i\in\{1,\ld,r\}$,
then $\bigl|\int_{Y}f\,d\mu_{Y}-\int_{X}f\,d\mu_{X}\bigr|<\eps$
for all functions $f$ on $X$ with $\sup_{x\neq y}|f(x)-f(y)|/\dist(x,y)\leq C$.
\endproposition
\npar
(This is in complete analogy with the situation on tori:
if $X$ is a torus,
for any $\eps>0$ there are only finitely many subtori $V_{1},\ld,V_{r}$
such that any subtorus $Y$ of $X$ that contains $0$
and is not contained in $\bigcup_{i=1}^{r}V_{i}$
is $\eps$-dense and ``$\eps$-uniformly distributed'' in $X$.)
\section{2}{Nilsequences, null-sequences, and generalized polynomials}

We will deal with the space $\Li=l^{\infty}(\Z^{d})$ 
of bounded sequences $\phi\col\Z^{d}\ra\C$, 
with the norm $\|\phi\|=\sup_{n\in\Z^{d}}|\phi(n)|$.

{\it A basic $r$-step nilsequence\/} is an element of $\Li$ 
of the form $\psi(n)=f(\eta(n)x)$, $n\in\Z^{d}$, 
where $x$ is a point of an $r$-step nilmanifold $X=G/\Gam$,
$\eta$ is a homomorphism $\Z^{d}\ra G$, and $f\in C(X)$.
We will denote the algebra of basic $r$-step nilsequences by $\cNb_{r}$, 
and the algebra $\bigcup_{r\in\N}\cNb_{r}$ of all basic nilsequences by $\cNb$.
We will denote the closure of $\cNb_{r}$, $r\in\N$, in $\Li$ by $\cN_{r}$,
and the closure of $\cNb$ by $\cN$;
the elements of these algebras will be called {\it $r$-step nilsequences\/}
and, respectively, {\it nilsequences}.

Given a polynomial sequence 
$g(n)=a_{1}^{p_{1}(n)}\kern-3mm\ld a_{k}^{p_{k}(n)}$, $n\in\Z^{d}$,
in a nilpotent group with $\deg p_{i}\leq s$ for all $i$,
we will say that $g$ has {\it naive degree\/} $\leq s$.
(The term ``degree'' was already reserved 
for another parameter of a polynomial sequence.)
We will call a sequence of the form $\psi(n)=f(g(n)x)$,
where $x$ is a point of an $r$-step nilmanifold $X=G/\Gam$,
$g$ is a polynomial sequence of naive degree $\leq s$ in $G$, 
and $f\in C(X)$,
{\it a basic polynomial $r$-step nilsequence of degree $\leq s$}.
We will denote the algebra of basic polynomial $r$-step nilsequences of degree $\leq s$
by $\cNb_{r,s}$ and the closure of this algebra in $\Li$ by $\cN_{r,s}$.
It is shown in \brfr{L2} (see proof of Theorem~B$^*$)
that any basic polynomial $r$-step nilsequence of degree $\leq s$
is a basic $l$-step nilsequence, where $l=2rs$;
we introduce this notion here only in order to keep trace 
of the parameters $r$, $s$ of the ``origination'' of a nilsequence.
So, for any $r$ and $s$, $\cNb_{r,s}\sle\cNb_{2rs}$;
since also $\cNb_{r}\sle\cNb_{r,1}$, we have $\bigcup_{r,s\in\N}\cNb_{r,s}=\cNb$.

We will also need the following lemma,
saying, informally, that the operation of ``alternation'' of sequences 
preserves the algebras of nilsequences:
\lemma{P-alter}{2.1}{}
Let $k\in\N$
and let $\psi_{i}\in\cNb_{r,s}$ (respectively, $\cNb_{r}$), $i\in\{0,\ld,k-1\}^{d}$.
Then the sequence $\psi$ defined by
$\psi(n)=\psi_{i}(m)$ for $n=km+i$ with $m\in\Z^{d}$, $i\in\{0,\ld,k-1\}^{d}$,
is also in $\cNb_{r,s}$ (respectively, $\cNb_{r})$).
\endlemma
\proof{}
Put $I=\{0,\ld,k-1\}^{d}$.
For each $i\in I$,
let $X_{i}=G_{i}/\Gam_{i}$ be the $r$-step nilmanifold,
$g_{i}$ be the polynomial sequence in $G_{i}$ of naive degree $\leq s$,
$x_{i}\in X_{i}$ be the point,
and $f_{i}\in C(X_{i})$ be the function
such that $\psi_{i}(n)=f(g_{i}(n)x_{i})$, $n\in\Z^{d}$.
If, for some $i$, $G_{i}$ is not connected,
it is a factor-group of a free $r$-step nilpotent group
with both continuous and discrete generators,
which, in its turn, is a subgroup of a free $r$-step nilpotent group 
with only continuous generators;
thus after replacing, if needed, $X_{i}$ by a larger nilmanifold
and extending $f_{i}$ to a continuous function on this nilmanifold
we may assume that every $G_{i}$ is connected and simply-connected.
In this case for any element $b\in G_{i}$ and any $l\in\N$
an $l$-th root $b^{1/l}$ exists in $G_{i}$,
and thus the polynomial sequence $b^{p(n)}$ in $G_{i}$ makes sense
even if a polynomial $p$ has non-integer rational coefficients.
Thus, for each $i\in I$
we may construct a polynomial sequence $g'_{i}$ in $G_{i}$,
of the same naive degree as $g_{i}$,
such that $g'_{i}(km+i)=g_{i}(m)$ for all $m\in\Z^{d}$.
Put $G=\Z^{d}\times\prod_{i\in I}G_{i}$,
$X=(\Z/(k\Z))^{d}\times\prod_{i\in I}X_{i}$,
$g(n)=\bigl(n,(g'_{i}(n),\ i\in I)\bigr)$ for $n\in\Z^{d}$,
$x=\bigl(0,(x_{i},\ i\in I)\bigr)\in X$,
and $f\bigl(i,(y_{i},\ i\in I)\bigr)=f_{i}(y_{i})$ for $(i,(y_{i},\ i\in I))\in X$.
Then $X$ is an $r$-step nilmanifold, $f\in C(X)$,
and thus the sequence
$\psi(n)=f(g(n)x)=f_{i}(g'_{i}(n)x_{i})=f_{i}(g_{i}(m)x_{i})=\psi_{i}(m)$,
$n=km+i$, $m\in\Z^{d}$, $i\in I$,
is in $\cN_{r,s}$.
\endproof

A set $S\sln\Z^{d}$ is said to be of {\it uniform {\rm(or {\it Banach})} density zero}
if for any F{\o}lner sequence $(\Phi_{N})_{N=1}^{\infty}$ in $\Z^{d}$, 
$\lim_{N\ras\infty}|S\cap\Phi_{N}|/|\Phi_{N}|=0$. 
A sequence $(\om_{n})_{n\in\Z^{d}}$ in a topological space $\Om$
converges to $\om\in\Om$ {\it in uniform density\/}
if for every neighborhood $U$ of $\om$
the set $\bigl(\{n\in\Z^{d}:\om_{n}\not\in U\}\bigr)$ is of uniform density zero.

We will say that a sequence $\lam\in\Li$ is {\it a null-sequence\/}
if it tends to zero in uniform density.
This is equivalent to 
$\lim_{N\ras\infty}\frac{1}{|\Phi_{N}|}\sum_{n\in\Phi_{N}}|\lam(n)|=0$
for any F{\o}lner sequence $(\Phi_{N})_{N=1}^{\infty}$ in $\Z^{d}$,
and is also equivalent to
$\lim_{N\ras\infty}\frac{1}{|\Phi_{N}|}\sum_{n\in\Phi_{N}}|\lam(n)|^{2}=0$
for any F{\o}lner sequence $(\Phi_{N})_{N=1}^{\infty}$ in $\Z^{d}$.
We will denote the set of (bounded) null-sequences by $\cZ$.
$\cZ$ is a closed ideal in $\Li$.

The algebra $\cZ$ is orthogonal to the algebra $\cN$ in the following sense: 
\lemma{P-orto}{2.2}{}
For any $\psi\in\cN$ and $\lam\in\cZ$,
$\|\psi+\lam\|\geq\|\psi\|$.
\endlemma
\proof{}
Let $c\geq\|\psi+\lam\|$.
Nilsystems are distal systems (see, for example, \brfr{L1}),
which implies that every point returns to any its neighborhood regularly.
It follows that if $|\psi(n)|>c$ for some $n$,
then the set $\bigl\{n\in\Z^{d}:|\psi(n)>c|\bigr\}$ has positive lower density,
and then $\psi(n)+\lam(n)>c$ for some $n$, contradiction.
Hence, $|\psi(n)|\leq c$ for all $n$.
\endproof
It follows that $\cN\cap\cZ=0$.

We will denote the algebras $\cNb_{r}+\cZ$, $\cN_{r}+\cZ$,
$\cNb_{r,s}+\cZ$, $\cN_{r,s}+\cZ$, $\cNb+\cZ$, and $\cN+\cZ$
by $\cMb_{r}$, $\cM_{r}$, $\cMb_{r,s}$, $\cM_{r,s}$, $\cMb$, and $\cM$ respectively.

\rfr{Lemma}{2.2} implies:
\lemma{P-clos}{2.3}{}
The algebras $\cM$, $\cM_{r}$, and $\cM_{r,s}$, $r,s\in\N$, are all closed,
and the projections $\cM\ra\cN$, $\cM\ra\cZ$ are continuous.
\endlemma
\proof{}
If a sequence $(\psi_{n}+\lam_{n})$ with $\psi_{n}\in\cN$, $\lam_{n}\in\cZ$, 
converges to $\phi\in\Li$,
then since $\|\psi_{n}\|\leq\|\psi_{n}+\lam_{n}\|$ for all $n$,
the sequence $(\psi_{n})$ is Cauchy, and so converges to some $\psi\in\cN$.
Then $(\lam_{n})$ also converges, to some $\lam\in\cZ$,
and so $\phi=\psi+\lam\in\cM$.
All the same is true for $\cM_{r}$ and $\cM_{r,s}$, instead of $\cM$, for all $r$ and $s$.

For $\phi_{1}=\psi_{1}+\lam_{1}$ and $\phi_{2}=\psi_{2}+\lam_{2}$
with $\psi_{1},\psi_{2}\in\cN$ and $\lam_{1},\lam_{2}\in\cZ$
we have $\|\psi_{1}-\psi_{2}\|\leq\|\phi_{1}-\phi_{2}\|$,
so the projection $\cM\ra\cN$, $\psi+\lam\mapsto\psi$, is continuous,
and so the projection $\cM\ra\cZ$, $\psi+\lam\mapsto\lam$, is also continuous.
\endproof

{\it Generalized polynomials\/} on $\Z^{d}$ 
are the functions on $\Z^{d}$ constructed from ordinary polynomials
using the operations of addition, multiplication, 
and the operation of taking the integer part.
A function $h$ on a nilmanifold $X$ is said to be {\it piecewise polynomial\/}
if it can be represented in the form
$h(x)=q_{i}(x)$, $x\in Q_{i}$, $i=1,\ld,k$, 
where $X=\bigcup_{i=1}^{k}Q_{i}$ is a finite partition of $X$
and, in Malcev coordinates on $X$, 
for every $i$ the set $Q_{i}$ is defined by a system of polynomial inequalities
and $q_{i}$ is a polynomial function.
(Since multiplication in a nilpotent Lie group is polynomial,
this definition does not depend on the choice of coordinates on $X$; see \brfr{BL}.)
It was shown in \brfr{BL} (and also, in a simpler way, in \brfr{L7}),
that a sequence $\ups\in\Li$ is a generalized polynomial
iff there is a nilmanifold $X=G/\Gam$, a piecewise polynomial function $h$ on $X$,
a polynomial sequence $g$ in $G$, and a point $x\in X$
such that $\ups(n)=h(g(n)x)$, $n\in\Z^{d}$.

Let $\cPb$ be the algebra of bounded generalized polynomials on $\Z^{d}$
and $\cP$ be the closure of $\cPb$ in $\Li$.
Since (by the Weierstrass approximation theorem)
any continuous function on a compact nilmanifold $X$
is uniformly approximable by piecewise polynomial functions,
any basic nilsequence is uniformly approximable by bounded generalized polynomials,
and so, is contained in $\cP$.
Hence, $\cN\sln\cP$, and $\cM\sln\cP+\cZ$.
The inverse inclusion does not hold,
since not all piecewise polynomial functions are uniformly approximable by continuous functions;
however, they are -- on the complement of a set of arbitrarily small measure,
which implies that generalized polynomials are also approximable by nilsequences, --
in a certain weaker topology in $\Li$.
\section{3}{Distribution of a polynomial sequence of subnilmanifolds}

Let $Y$ be a connected subnilmanifold of the (connected) nilmanifold $X$,
and let $g(n)$, $n\in\Z^{d}$, be a polynomials sequence in $G$.
We will investigate how the sequence $g(n)Y$ of sunilmanifolds of $X$
is distributed in $X$.

\proposition{P-YtoZ}{3.1}{}
Let $X=G/\Gam$ be a connected nilmanifold,
let $Y$ be a connected subnilmanifold of\/ $X$,
and let $g\col\Z^{d}\ra G$ be a polynomial sequence in $G$ with $g(0)=\ed{G}$.
Assume that $g(\Z^{d})Y$ is dense in $X$,
and that $G$ is generated by $\Gc$ and the range $g(\Z^{d})$ of $g$.
Let $Z$ be the normal closure of\/ $Y$ in $X$;
then for any $f\in C(X)$,
$\lam(n)=\int_{g(n)Y}f\,d\mu_{g(n)Y}-\int_{g(n)Z}f\,d\mu_{g(n)Z}$, $n\in\Z^{d}$, is a null-sequence.
\endproposition
\proof{}
Let $f\in C(X)$ and let $\eps>0$;
we have to show that the set
$\bigl\{n\in\Z^{d}:\bigl|\int_{g(n)Y}f\,d\mu_{g(n)Y}-\int_{g(n)Z}f\,d\mu_{g(n)Z}\bigr|\geq\eps\bigr\}$
has zero uniform density in $\Z^{d}$.
After replacing $f$ by a close function we may assume that $f$ is Lipschitz,
so that $C=\sup_{x\neq y}|f(x)-f(y)|/\dist(x,y)$ is finite.
Choose Malcev's coordinates in $\Gc$,
and let $Q\sln\Gc$ be the corresponding fundamental cube.
Since $Z$ is normal in $X$, 
$aZ=bZ$ whenever $a=b\mod\Gam$,
and $\bigcup_{a\in Q}aZ$ is a partition of $X$.

We first want to determine for which $a\in G$
one has $\bigl|\int_{aY}f_{aY}\,d\mu_{aY}-\int_{aZ}f\,d\mu_{aZ}\bigr|\geq\eps$.
For every $b\in Q$,
by \rfr{Proposition}{1.1}, applied to the nilmanifold $bZ$,
there exist proper subnilmanifolds $V_{b,1},\ld,V_{b,r_{b}}$ of $Z$
such that $\bigl|\int_{W}f\,d\mu_{W}-\int_{bZ}f_{bZ}\,d\mu_{bZ}\bigr|<\eps/2$
whenever $W$ is a subnilmanifold of $bZ$ 
with $b\in W\not\sle bV_{b,i}$, $i=1,\ld,r_{b}$.
By continuity, for each $b\in Q$
there exists a neighborhood $U_{b}$ of $b$ such that for all $a\in U_{b}$,
$\bigl|\int_{W}f\,d\mu_{W}-\int_{aZ}f\,d\mu_{aZ}\bigr|<\eps$
whenever $a\in W\sle aZ$ and $W\not\sle aV_{b,i}$, $i=1,\ld,r_{b}$.
Using the compactness of the closure of $Q$,
we can choose $b_{1},\ld,b_{l}\in Q$ such that $Q\sle\bigcup_{j=1}^{l}U_{b_{j}}$;
let $V=\bigcup_{\sdup{j=1,\ld,l}{i=1,\ld,r_{j}}}V_{b_{j},i}$.
Then for any $b\in Q$,
for any subnilmanifold $W$ of $bZ$ with $b\in W\not\sle V$
one has $\bigl|\int_{W}f\,d\mu_{W}-\int_{bZ}f\,d\mu_{bZ}\bigr|<\eps$.
Now let $a\in G$, and let $b\in Q$ be such that $a=b\mod\Gam$.
Then $aY\sle aZ=bZ$,
thus, if $aY\not\sle bV$,
then $\bigl|\int_{aY}f\,d\mu_{aY}-\int_{aZ}f\,d\mu_{aZ}\bigr|<\eps$.
Hence, $\bigl|\int_{aY}f\,d\mu_{aY}-\int_{aZ}f\,d\mu_{aZ}\bigr|\geq\eps$
only if $(a1_{X},aY)\sle(b1_{X},bV)$ for some $b\in Q$.

Let $N=\bigl\{(b1_{X},bV),\ b\in Q\bigr\}$;
we have to prove that the the set $\bigl\{n\in\Z^{d}:(g(n)1_{X},g(n)Y)\dsc\sle N\bigr\}$
has zero uniform density in $\Z^{d}$.
For this purpose we are going to find the closure
of the sequence $\tY_{n}=(g(n)1_{X},g(n)Y)$, $n\in\Z^{d}$,
(more exactly, of the union $\bigcup_{n\in\Z^{d}}\tY_{n}$),
the orbit of the subnilmanifold $\tY=(1_{X},Y)$ of $X\times X$
under the polynomial action $(g(n),g(n))$, $n\in\Z^{d}$.
Assume for simplicity 
that the closure $R$ of the orbit $\bigl\{g(n)1_{X},\ n\in\Z^{d}\bigr\}$ is connected,
and let $P$ be the closed connected subgroup of $G$ such that $\pi(P)=R$.
(If $R$ is disconnected 
we pass to a sublattice of $\Z^{d}$ and its cosets
to deal with individual connected components of $R$.)
We will also assume that $Y\ni 1_{X}$.
\lemma{P-orb}{3.2}{}
The closure of the sequence $(\tY_{n})$ is the subnilmanifold
$D=\bigl\{(a1_{X},aZ),\ a\in P\bigr\}=\bigl\{(a1_{X},aZ),\ a\in P\cap Q\bigr\}$ 
of $X\times X$.
\endlemma
\proof{}
Let $L$ be the closed connected subgroup of $G$ such that $\pi(L)=Z$,
and let $K=\bigl\{(a,au),\ a\in P,\ u\in L\bigr\}$;
since $L$ is normal in $G$, $K$ is a (closed rational) subgroup of $G\times G$,
and we have $D=K/((\Gam\times\Gam)\cap K)$.

For any $n\in\Z^{d}$ we have $\tY_{n}\sle D$
(since $g(n)1_{X}\in R$, so $g(n)\in P\Gam$,
so $g(n)L\sle P\Gam L=PL\Gam$),
and we have to show that the sequence $(\tY_{n})$ is dense in $D$.
For this it suffices to prove that the image of this sequence is dense
in the maximal torus $T=\lfa{D}{[K,K]}$ of $D$.
Since $L$ is normal,
we have $[K,K]=\bigl\{(a,au),\ a\in[P,P],\ u\in[P,L][L,L]\bigr\}$,
and the torus $T$ is the factor of the commutative group $K/[K,K]$
by the image $\Lam$ in this group of the lattice $\Gam\times\Gam$.
Let $H$ be the closed connected subgroup of $G$ such that $\pi(H)=Y$;
then $L=H[H,G]$, so
$$
K/[K,K]=\bigl\{(a,avw),\ a\in P,\ v\in H,\ w\in[H,G]\bigr\}/
\bigl\{(a,au),\ a\in[P,P],\ u\in[P,L][L,L]\bigr\}.
$$

By assumption, $G$ is generated by $\Gc$ and $g$.
Since the orbit $\{g(n)Z,\ n\in\Z^{d}\}$ is dense in $X$ and $Z$ is normal,
the orbit $\{g(n)1_{X/Z},\ n\in\Z^{d}\}$ is dense in $X/Z$,
so $P/(P\cap L)=\Gc/L$, so $\Gc=PL$.
Hence, $[H,G]=[H,g][H,P][H,L]$.
For any $n$, $g(n)=u_{n}\gam_{n}$ for some $u_{n}\in P$ and $\gam_{n}\in\Gam$, 
thus, modulo $[P,L][L,L]$, 
the group $[H,G]$ is generated by $\{[H,\gam_{n}],\ n\in\Z^{d}\}$.

The closure $B$ of the image of the sequence $(\tY_{n})$ in $T$ is a subtorus of $T$.
Since the sequence $g(n)1_{X}$ is dense in $R=\pi(P)$,
the subtorus $T_{1}=\{(a,a),\ a\in P\}/([K,K]\Lam)$ of $T$
is the closure of the image of the sequence $(g(n)1_{X},g(n)1_{X})$
and so, is contained in $B$.
Also, the subtorus $T_{2}=\{(1_{X},u),\ u\in H\}/([K,K]\Lam)$ is contained in $B$.
Finally, for $n\in\Z^{d}$ and $c\in H$ we have
$$
(g(n),g(n)c)=(u_{n},u_{n}\gam_{n}c)=(u_{n}1_{X},u_{n}c[c,\gam_{n}^{-1}]\gam_{n}).
$$
Taken modulo $[K,K]\Lam$, these elements of $B$ generate $T$ modulo $T_{1}+T_{2}$,
so $B=T$.
\endproof

It follows that the sequence $\tY_{n}$, $n\in\Z^{d}$, is well distributed in $D$.
The set $N=\bigl\{(b1_{X},bV),\ b\in Q\bigr\}$
is a compact subset of $D$ of zero measure,
thus, the set $\bigl\{n\in\Z^{d}:(g(n)1_{X},g(n)Y)\sle N\bigr\}$
has zero uniform density in $\Z^{d}$.
\endproof

\theorem{P-Ynils}{3.3}{}
Let $X=G/\Gam$ be an $r$-step nilmanifold,
let $Y$ be a subnilmanifold of $X$,
let $g$ be a polynomial sequence in $G$ of naive degree $\leq s$,
let $f\in C(X)$.
Then the sequence $\phi(n)=\int_{g(n)Y}f\,d\mu_{g(n)Y}$, $n\in\Z^{d}$,
is contained in $\cMb_{r,s}$.
\endtheorem
\proof{}
We may assume that $Y\ni\ed{X}$.
After replacing $f$ by $f(g(0)x)$, $x\in X$,
we may assume that $g(0)=1_{X}$.
We may also replace $X$ by the closure of the orbit $g(\Z^{d})Y$,
and we may assume that $G$ is generated by $\Gc$ and the range of $g$.

First, let $X$ and $Y$ be both connected.
Let $Z$ be the normal closure of $Y$ in $X$;
then by \rfr{Proposition}{3.1}, 
$\phi(n)=\int_{g(n)Z}f\,d\mu_{g(n)Z}+\lam_{n}$, $n\in\Z^{d}$,
with $\lam\in\cZ$. 
Define $\hX=X/Z$, $\hx=\{Z\}\in\hX$, and $\hf=E(f|\hX)\in C(\hX)$;
then $\int_{g(n)Z}f\,d\mu_{g(n)Z}=\hf(g(n)\hx)$, $n\in\Z^{d}$,
and the sequence $\hf(g(n)\hx)$, $n\in\Z^{d}$, is in $\cNb_{r,s}$,
so $\phi\in\cMb_{r,s}$.

Now assume that $Y$ is connected but $X$ is not.
Then, by \brfr{L2},
there exists $k\in\N$ such that $\ovr{g(k\Z^{d}+i)Y}$ is connected
for every $i\in\{0,\ld,k-1\}^{d}$.
Thus, for every $i\in\{0,\ld,k-1\}^{d}$, $\phi(kn+i)\in\cMb_{r,s}$,
and the assertion follows from \rfr{Lemma}{2.1}.

Finally, if $Y$ is disconnected 
and $Y_{1},\ld,Y_{r}$ are the connected components of $Y$,
then $\int_{g(n)Y}f\,d\mu_{g(n)Y}
=\sum_{j=1}^{r}\int_{g(n)Y_{j}}f\,d\mu_{g(n)Y_{j}}$, $n\in\Z^{d}$,
and the result holds since it holds for $Y_{1},\ld,Y_{r}$.
\endproof
\section{4}{Integrals of null- and of nil-sequences}

On $\Li$ and, thus, on $\cN$, $\cZ$ and $\cM$ 
we will assume the Borel $\sig$-algebra induced by the weak topology.

We start with integration of null-sequences:
\lemma{P-IntZ}{4.1}{}
Let $(\Om,\nu)$ be a measurable space
and let $\Om\ra\cZ$, $\om\mapsto\lam_{\om}$, be an integrable mapping.
Then the sequence $\lam(n)=\int_{\Om}\lam_{\om}(n)\,d\nu$ is in $\cZ$ as well.
\endlemma
\npar
(We say that a mapping $\Psi\col\Om\ra\Li$ is {\it integrable\/}
if it is measurable and $\int_{\Om}\|\Psi\|\,d\nu<\infty$.)
\proof{}
For each $\om\in\Om$,
$\lim_{N\ras\infty}\frac{1}{|\Phi_{N}|}\sum_{n\in\Phi_{N}}|\lam_{\om}(n)|=0$
for any F{\o}lner sequence $(\Phi_{N})_{N=1}^{\infty}$ in $\Z^{d}$.
By the dominated convergence theorem,
\lequ{
\limsup_{N\ras\infty}\frac{1}{|\Phi_{N}|}\sum_{n\in\Phi_{N}}|\lam(n)|
=\limsup_{N\ras\infty}\frac{1}{|\Phi_{N}|}\sum_{n\in\Phi_{N}}\Bigl|\int_{\Om}\lam_{\om}(n)\,d\nu\Bigr|
\-\\\-
\leq\lim_{N\ras\infty}\int_{\Om}\frac{1}{|\Phi_{N}|}\sum_{n\in\Phi_{N}}|\lam_{\om}(n)|\,d\nu
=\int_{\Om}\lim_{N\ras\infty}\frac{1}{|\Phi_{N}|}\sum_{n\in\Phi_{N}}|\lam_{\om}(n)|\,d\nu
=0.
}
So, $\lam\in\cZ$.
\endproof

For nilsequences we have:
\proposition{P-IntN}{4.2}{}
Let $(\Om,\nu)$ be a measure space
and let $\Om\ra\cN$, $\om\mapsto\phi_{\om}$, be an integrable mapping.
Then the sequence $\phi(n)=\int_{\Om}\phi_{\om}(n)\,d\nu$ belongs to $\cM$.
(If, for some $r$, $\phi_{\om}\in\cN_{r}$ for all $\om$, then $\phi\in\cM_{r}$.)
\endproposition

To simplify notation, let us start with the case $d=1$.
We are going to reduce \rfr{Proposition}{4.2} 
to a statement concerning a sequence of measures on a nilmanifold.%
\footnote{$^{1}$}{The argument that follows has been changed;
I thank B.~Host for pointing to me out a mistake in the previos version of the paper.}
Since $\cM$ is closed in $\Li$,
we are allowed to replace the mapping $\phi_{\om}$ from $\Om$ to $\cN$
by a close mapping $\phi'_{\om}$:
we are done if for any $\eps>0$ we can find a mapping $\Om\ra\cN$, $\om\mapsto\phi'_{\om}$,
with $\|\int_{\Om}\|\phi_{\om}-\phi'_{\om}\|_{\Li}\,d\nu<\eps$
and such that the assertion of \rfr{Proposition}{4.2} holds for $\phi'_{\om}$.
Fix $\eps>0$.
First, after replacing $\Om$ by $\Om'\sle\Om$ with $\nu(\Om')<\infty$
such that $\int_{\Om\sm\Om'}\|\phi_{\om}\|_{\Li}\,d\nu<\eps$, 
we may assume that $\nu(\Om)<\infty$.
Next, since the set $\cNb$ of basic nilsequences is dense in $\cN$,
we may replace the nilsequences $\phi_{\om}$ by close basic nilsequences,
if we manage to do this in a measurable way.
We will, as we may, deal with $\R$-valued nilsequences.
Let $X=G/\Gam$ be a nilmanifold
where $G$ is a simply connected nilpotent Lie group and $\Gam$ is a lattice in $G$,
and let $\pi\col G\ra X$ be the projection.
We may assume that $G$ has the same number of connected components as $X$,
then $G$ is homeomorphic to $\R^{\dim G}\times F$, where $F$ is a finite set,
with $\Gam$ corresponding to $\Z^{\dim G}$;
this homeomorphism induces a natural metric on $G$ and on $X$.
For $k\in\N$
let $Q_{k}$ be the set of elements of $G$ at the distance $\leq k$ from $1_{G}$
and let $L_{k}$ be the set of Lipschitz functions on $X$
with Lipshitz constant $k$ and of modulus $\leq k$.
The subset $Q_{k}\times L_{k}$ of $G\times C(X)$ is compact;
the ``nilsequence reading'' mapping $\Psi\col G_{k}\times C(X_{k})\ra\cN$,
$\Psi(a,h)(n)=h(\pi(a^{n}))$,
is continuous with respect to the weak topology on $\cN$;
thus the set $\cL_{X,k}=\Psi(Q_{k}\times L_{k})\sln\cNb$ is compact in this topology.
Fix a countable set $S$ dense in $\cL_{X,k}$ in the weak topology and enumerate it.
Let $\phi\in\cN$.
For each $j\in\N$ let $\psi_{j}$ be the element of $S$ for which\\
(i) the sum $\sum_{n=-j}^{j}|\phi(n)-\psi_{j}(n)|$ is minimal;\\
(ii) among the elements of $S$ for which (i) holds, 
the vector $\bigl(\psi(0),\psi(-1),\psi(1),\ld,\psi(-j),\psi(j)\bigr)$
is minimal for $\psi=\psi_{j}$ with respect to the lexicographic order;\\
(iii) and among the elements of $S$ for which (i) and (ii) hold, 
$\psi_{j}$ has the minimal number under the ordering of $S$.\\
Put $\zeta_{X,k,j}(\phi)=\psi_{j}$;
then $\zeta_{X,k,j}$ is a measurable mapping $\cN\ra\cL_{X,k}$.
For any $\phi\in\cN$ the sequence $\psi_{j}=\zeta_{X,k,j}(\phi)$ converges in $\cL_{k}$:
indeed, $\cL_{X,k}$ is compact,
and any convergent subsequence of this sequence converges to the same element of $\cL_{X,k}$,
namely, to $\psi\in\cQ_{k}$ which is closest to $\phi$ in the $\Li$-norm,
and among such, which is minimal with respect to the lexicographic order of its entries.
Put $\zeta_{X,k}(\phi)=\lim_{j\ras\infty}\zeta_{X,k,j}(\phi)$, $\phi\in\cN$;
then $\zeta_{X,k}$ is a measurable mapping $\cN\ra\cL_{X,k}$
that maps each nilsequence to a closest in $\Li$-norm element of $\cL_{X,k}$.
It also follows that the function $\dX_{X,k}(\phi)=\min_{\psi\in\cL_{X,k}}\|\phi-\psi\|_{\Li}$
is measurable on $\cN$.

In each class of isomorphic nilmanifolds choose a representative $X$
(along with $G$, $\Gam$, a homeomorphism $G\ra\R^{\dim G}\times F$, 
and a metric on $G$ and $X$);
let $\cX$ be the set of these representatives.
Since there exists only countably many nonisomorphic nilmanifolds, $\cX$ is countable.
Introduce a well ordering of $\cX$ satisfying $X'<X$ when $\dim X'<\dim X$.
For every $X\in\cX$
put $\Om_{X,k}=\bigl\{\om\in\Om:\dX_{X,k}(\phi_{\om})<\eps/\nu(\Om)\bigr\}$
and $\Om_{X}=\bigcup_{k=1}^{\infty}\Om_{X,k}$;
these are measurable subsets of $\Om$.
The union $\bigcup_{X\in\cX}\bigcup_{k=1}^{\infty}\cL_{X,k}$ is dense in $\cN$,
thus $\bigcup_{X\in\cX}\bigcup_{k=1}^{\infty}\Om_{X,k}=\Om$.
Next define $\Om'_{X}=\Om_{X}\sm\bigcup_{X'<X}\Om_{X'}$, $X\in\cX$;
these are disjoint sets that partition $\Om$.
Finally, for each $X\in\cX$ and $k\in\N$ put $\Om'_{X,k}=\Om'_{X}\sm\bigcup_{k'<k}\Om'_{X,k'}$.
Now, for $\om\in\Om$ define $\psi_{\om}=\zeta_{X,k}(\phi_{\om})$ 
when $\om\in\Om'_{X,k}$, $X\in\cX$, $k\in\N$;
then $\om\mapsto\psi_{\om}$ is a measurable mapping $\Om\ra\cNb$
with $\|\psi_{\om}-\phi_{\om}\|_{\Li}<\eps/\nu(\Om)$ for all $\om\in\Om$.
We may now replace $\phi_{\om}$ by $\psi_{\om}$, $\om\in\Om$;
moreover, we may also deal with the sets $\Om'_{X}$ separately,
and therefore assume that $\phi_{\om}$, $\om\in\Om$, 
are all read off the same nilmanifod $X=G/\Gam$:
$\phi_{\om}=\Psi(a,h)$ with $a\in G$ and $h$ being a Lipschitz function on $X$.
(And, in addition, by our construction, 
$\phi_{\om}$ is not readable off any nilmanifold $X'$ with $X'<X$.)

We now claim that for each $\om\in\Om$, 
$\phi_{\om}$ has only countably many preimages under this mapping;
by Lusin's theorem about the existence of a measurable section, 
this will imply that we can choose elements $a_{\om}\in Q_{k}$, $h_{\om}\in L_{k}$ 
with $\phi_{\om}=\Psi(a_{\om},h_{\om})$, $\om\in\Om$,
such that the mapping $\om\mapsto(a_{\om},h_{\om})$ is measurable.
We get use of the following fact:
\lemma{P-orfac}{4.2a}{}
Let nilmanifolds $X_{i}=G_{i}/\Gam_{i}$, elements $a_{i}\in G_{i}$, 
and functions $h_{i}\in C(X_{i})$, $i=1,2$,
be such that the orbits $\{\pi_{i}(a_{i}^{n})\}_{n\in\Z}$,
where $\pi_{i}$ are the projections $G_{i}\ra X_{i}$, are dense in $X_{i}$, $i=1,2$,
and the triples $(X_{1},a_{1},h_{1})$ and $(X_{2},a_{2},h_{2})$
produce the same nilsequence: 
$\phi(n)=h_{1}(\pi_{1}(a_{1}^{n}))=h_{2}(\pi_{2}(a_{2}^{n}))$, $n\in\Z$.
Then there exists a common factor $(\hX,\ha,\hh)$ of $(X_{1},a_{1},h_{1})$ and $(X_{2},a_{2},h_{2})$
such that $\phi(n)=\hh(\hpi(\ha^{n}))$ (where $\hpi$ is the projection $\hG\ra\hX$).
\endlemma
\proof{}
Let $\tG=G_{1}\times G_{2}$, $\tX=X_{1}\times X_{2}$, $\tpi=\pi_{1}\times\pi_{2}\col\tG\ra\tX$.
Let $\ta=(a_{1},a_{2})\in\tG$,
and let $Y$ be the closure of the orbit of $\ta$ in $\tX$, $Y=\ovr{\{\tpi(\ta^{n}),\ n\in\Z\}}$.
Let $p_{1}$ and $p_{2}$ be the projections of $Y$ to $X_{1}$ and to $X_{2}$ respectively.
Let $Y_{1}\sle X_{1}$ and $Y_{2}\sle X_{2}$ be ``the fibers'' 
of the projections $p_{2}$ and $p_{1}$ respectively:
$p_{2}^{-1}(1_{X_{2}})=Y_{1}\times\{1_{X_{2}}\}$
and $p_{1}^{-1}(1_{X_{1}})=Y_{2}\times\{1_{X_{1}}\}$.
For any $n\in\Z$ we have $\phi(n)=h_{1}(p_{1}(\tpi(\ta^{n})))=h_{2}(p_{2}(\tpi(\ta^{n})))$;
since the orbit $\{\tpi(\ta^{n})\}_{n\in\Z}$ is dense in $Y$,
this implies that $h_{1}\comp p_{1}=h_{2}\comp p_{2}$,
so $h_{1}$ is constant on the fibers $b_{1}Y_{1}$, $b_{1}\in G_{1}$, of $p_{2}$
and $h_{2}$ is constant on the fibers $b_{2}Y_{2}$, $b_{2}\in G_{2}$, of $p_{1}$,
and therefore we may factorize $X_{1}$ by $Y_{1}$ and $X_{2}$ by $Y_{2}$ 
(and $Y$ by $Y_{1}\times Y_{2}$) to get the same factor $(\hX,\ha,\hh)$.
\endproof
Now, assume that for some $\om\in\Om$, 
$\phi_{\om}=\Psi(a_{1},h_{1})=\Psi(a_{2},h_{2})$, 
$a_{1},a_{2}\in G$ and $h_{1},h_{2}$ are Lipschitz functions on $X$.
Let, by \rfr{Lemma}{4.2a},
$(\hX,\ha,\hh)$ be a common factor of $(X,a_{1},h_{1})$ and $(X,a_{2},h_{2})$
such that $\phi_{\om}=\Psi(\ha,\hh)$.
Since $\phi_{\om}$ cannot be read off a nilmanifold $\hX$ with $\dim\hX<\dim X$,
there must be $\dim\hX=\dim X$.
However, for any pair $(\hX,\ha)$, ``nilmanifold with a translation'',
there are only countably many (up to isomorphism) pairs $(X,a)$ 
extending $(\hX,\ha)$ and with $\dim X=\dim\hX$.

Thus, we arrive at the following situation:
we have a nilmanifold $X=G/\Gam$ 
and a measurable function $\Om\ra G\times C(X)$, $\om\mapsto(a_{\om},h_{\om})$,
such that for every $\om\in\Om$ one has $\phi_{\om}(n)=h_{\om}(\pi(a_{\om}^{n}))$, $n\in\Z$.
Let $H(\om,x)=h_{\om}(x)$, $\om\in\Om$, $x\in X$;
then $\phi(n)=\int_{\Om}H(\om,\pi(a_{\om}^{n}))d\nu(\om)$, $n\in\Z$.
Choose a basis $f_{1},f_{2},\ld$ in $C(X)$;
the function $H$ is representable in the form
$H(\om,x)=\sum_{i=1}^{\infty}\theta_{i}(\om)f_{i}(x)$,
where convergence is uniform with respect to $x$ for any $\om$;
we are done if we prove the assertion for the functions $\theta_{i}(\om)f_{i}(x)$ instead of $H$.
So, let $\theta\in L^{1}(\Om)$ and $f\in C(X)$;
we have to show that the sequence $\phi(n)=\int_{\Om}\theta(\om)f(\pi(a_{\om}^{n}))\,d\nu(\om)$
is in $\cM$.
We may also assume that $\theta\geq 0$.
Let $\tau\col\Om\ra G$ be the mapping defined by $\tau(\om)=a_{\om}$,
and let $\ro=\tau_{*}(\theta\nu)$;
then $\ro$ is a finite measure on $G$ and $\phi(n)=\int_{G}f(\pi(a^{n}))\,d\ro(a)$.
Thus, \rfr{Proposition}{4.2} will follow from the following:
\proposition{P-limMes}{4.3}{}
Let $X=G/\Gam$ be a nilmanifold,
let $\ro$ be a finite Borel measure on $G$,
and let $f\in C(X)$.
Then the sequence $\phi(n)=\int_{G}f(\pi(a^{n}))\,d\ro(a)$, $n\in\Z$, is in $\cMb$.
(If $X$ is an $r$-step nilmanifold, then $\phi\in\cMb_{r}$.)
\endproposition
\proof{}
We may and will assume that $X$ is connected.
Let $\tro=\pi_{*}(\ro)$;
we decompose $\tro$ in the following way:
\lemma{P-razr}{4.4}{}
There exists an at most countable collection $\cV$ of connected subnilmaniolds of $X$
(which may include $X$ itself and singletons)
and finite Borel measures $\ro_{V}$, $V\in\cV$, on $G$
such that $\ro=\sum_{V\in\cV}\ro_{V}$
and for every $V\in\cV$, $\supp(\tro_{V})\sle V$
and $\tro_{V}(W)=0$ for any proper subnilmanifold $W$ of $V$,
where $\tro_{V}=\pi_{*}(\ro_{V})$.
\endlemma
\proof{}
Let $\cV_{0}$ be the (at most countable) set of the singletons $V=\{x\}$ in $X$
(connected 0-dimensional subnilmanifolds of $X$)
for which $\tro(V)>0$.
For each $V\in\cV_{0}$ let $\ro_{V}$ be the restriction of $\ro$ to $\pi^{-1}(V)$
(that is, $\tro_{V}(A)=\ro(A\cap\pi^{-1}(V))$ for measurable subsets $A$ of $G$),
and let $\ro_{1}=\ro-\sum_{V\in\cV_{0}}\ro_{V}$ and $\tro_{1}=\pi_{*}(\ro_{1})$.
Now let $\cV_{1}$ be the (at most countable) set of connected 1-dimensional subnilmanifolds of $X$
for which $\tro_{1}(V)>0$,
for each $V\in\cV_{1}$ let $\ro_{V}$ be the restriction of $\ro_{1}$ to $\pi^{-1}(V)$,
and $\ro_{2}=\ro-\sum_{V\in\cV_{1}}\ro_{V}$, $\tro_{2}=\pi_{*}(\ro_{2})$.
(Note that for $V_{1},V_{2}\in\cV_{1}$, 
the subnilmanifold $V_{1}\cap V_{2}$, if nonempty, has dimension 0,
so $\tro_{1}(V_{1}\cap V_{2})=0$.)
And so on, by induction on the dimension of the subnilmanifolds;
at the end, we put $\cV=\bigcup_{i=0}^{\dim X}\cV_{i}$.
\endproof

By \rfr{Lemma}{2.3}, it suffices to prove the assertion for each of $\ro_{V}$ instead fo $\ro$.
So, we will assume that the measure $\ro$ is supported by a connected subnilmanifold $V$ of $X$
and $\ro(W)=0$ for any proper subnilmanifold $W$ of $V$.

First, let $V=X$:
\lemma{P-VeX}{4.5}{}
Let $\ro$ be a finite Borel measure on $G$
such that for $\tro=\pi_{*}(\ro)$
one has $\tro(W)=0$ for any proper subnilmanifold $W$ of $X$.
Then for any $f\in C(X)$
the sequence $\phi(n)=\int_{G}f(\pi(a^{n}))\,d\ro(a)$, $n\in\Z$, 
converges to $\int_{X}f\,d\mu_{X}$ in uniform density.
\endlemma
\proof{}
We may assume that $\int_{X}f\,d\mu_{X}=0$;
we then have to show that $\phi$ is a null-sequence.
Let $(\Phi_{N})$ be a  F{\o}lner sequence in $\Z$.
By the dominated convergence theorem we have
\equ{
\lim_{N\ras\infty}\frac{1}{|\Phi_{N}|}\sum_{n\in\Phi_{N}}|\phi(n)|^{2}
\ &=\lim_{N\ras\infty}\frac{1}{|\Phi_{N}|}\sum_{n\in\Phi_{N}}
\int_{G}f(\pi(a^{n}))\,d\ro(a)\int_{G}\bar{f}(\pi(b^{n}))\,d\ro(b)
\-\\
&=\lim_{N\ras\infty}\frac{1}{|\Phi_{N}|}\sum_{n\in\Phi_{N}}
\int_{G\times G}f(\pi(a^{n}))\bar{f}(\pi(b^{n}))\,d(\ro\times\ro)(a,b)
\-\\
&=\int_{G\times G}\lim_{N\ras\infty}\frac{1}{|\Phi_{N}|}\sum_{n\in\Phi_{N}}
f\otimes\bar{f}(\pi^{\times 2}(a^{n},b^{n}))\,d\ro^{\times2}(a,b)
\-\\
&=\int_{G\times G}F(a,b)\,d\ro^{\times2}(a,b),
\-}
where $\pi^{\times2}=\pi\times\pi$, $\ro^{\times2}=\ro\times\ro$,
and $F(a,b)=\lim_{N\ras\infty}\frac{1}{|\Phi_{N}|}\sum_{n\in\Phi_{N}}
f\otimes\bar{f}(\pi^{\times 2}(a^{n},b^{n}))$, $a,b\in G$.
For $a,b\in G$, 
if the sequence $u_{n}=\pi^{\times2}(a^{n},b^{n})$, $n\in\Z$, is well distributed in $X\times X$ 
then $F(a,b)=\int_{X\times X}f\otimes\bar{f}\,d\mu_{X\times X}
=\int_{X}f\,d\mu\int_{X}\bar{f}\,d\mu=0$.
So, $F(a,b)\neq 0$
only if the sequence $(u_{n})$ is not well distributed in $X\times X$,
which only happens if the point $\pi^{\times2}(a,b)$ is contained in a proper subnilmanifold $D$ 
of $X\times X$ with $1_{X\times X}\in D$.
So, 
$$
\lim_{N\ras\infty}\frac{1}{|\Phi_{N}|}\sum_{n\in\Phi_{N}}|\phi(n)|^{2}
\leq\sum_{D\in\cD}\int_{(\pi^{\times2})^{-1}(D)}|F(a,b)|\,d\ro^{\times2}(a,b),
$$
where $\cD$ is the (countable) set of proper subnilmanifolds of $X\times X$
containing $1_{X\times X}$.
Let $D\in\cD$;
then either for any $x\in X$ the fiber $W'_{x}=\{y\in X:(x,y)\in D\}$ of $D$ over $x$
is a proper subnilmanifold of $X$,
or for any $y\in X$ the fiber $W''_{y}=\{x\in X:(x,y)\in D\}$ of $D$ over $y$
is a proper subnilmanifold of $X$,
(or both).
Since, by our assumption, $\tro(W)=0$ for any proper subnilmanifold $W$ of $X$,
in either case $\tro^{\times2}(D)=0$,
so $\ro^{\times2}((\pi^{\times2})^{-1}(D))=0$.
Hence, $\lim_{N\ras\infty}\frac{1}{|\Phi_{N}|}\sum_{n\in\Phi_{N}}|\phi(n)|^{2}=0$,
which means that $\phi\in\cZ$.
\endproof
Thus, in this case, $\phi$ is a constant plus a null-sequence,
that is, $\phi\in\cMb$.

Let now $V$ be of the form $V=cY$, 
where $Y$ is a (proper) connected subnilmanifold of $X$ with $1_{X}\in Y$ and $c\in\Gc$.
We may and will assume that the orbit $\{c^{n}Y,\ n\in\Z\}$ of $Y$ is dense in $X$.
Let $Z$ be the normal closure of $Y$ in $X$.
In this situation the following generalization of \rfr{Lemma}{4.5} does the job:
\lemma{P-VnX}{4.6}{}
Let $Z$ be a normal subnilmanifold of $X$ 
and let $c\in G$ be such that $\{c^{n}Z, \ n\in\Z\}$ is dense in $X$.
Let $\ro$ be a finite Borel measure on $G$
such that for $\tro=\pi_{*}(\ro)$
one has $\supp(\tro)\sle cZ$
and $\tro(cW)=0$ for any proper normal subnilmanifold $W$ of $Z$.
Let $\phi(n)=\int_{G}f(\pi(a^{n}))\,d\ro(a)$, $n\in\Z$,
let $\hX=X/Z$, and let $\hf=E(f|\hX)$.
Then $\phi-\hf(\pi(c^{n}))\in\cZ$.
\endlemma
\proof{}
After replacing $f$ by $f-\hf$ we will assume that $E(f|\hX)=0$;
we then have to prove that $\phi$ is a null-sequence.
Let $L=\pi^{-1}(Z)$.
Let $(\Phi_{N})$ be a  F{\o}lner sequence in $\Z$;
as in \rfr{Lemma}{4.5}, we obtain
\equ{
\lim_{N\ras\infty}\frac{1}{|\Phi_{N}|}\sum_{n\in\Phi_{N}}|\phi(n)|^{2}
=\int_{G\times G}F(a,b)\,d\ro^{\times2}(a,b)
=\int_{(cL)\times(cL)}F(a,b)\,d\ro^{\times2}(a,b),
}
where $F(a,b)=\lim_{N\ras\infty}\frac{1}{|\Phi_{N}|}\sum_{n\in\Phi_{N}}
f\otimes\bar{f}(\pi^{\times 2}(a^{n},b^{n}))$, $a,b\in L$.
Let us ``shift'' $\ro$ to the origin,
by replacing it by $c^{-1}_{*}\ro(a)$, $a\in G$,
so that now $\supp(\ro)\sle L$, $\supp(\tro)\sle Z$, and 
\equ{
\lim_{N\ras\infty}\frac{1}{|\Phi_{N}|}\sum_{n\in\Phi_{N}}|\phi(n)|^{2}
=\int_{L\times L}F(a,b)\,d\ro^{\times2}(a,b),
}
where $F(a,b)=\lim_{N\ras\infty}\frac{1}{|\Phi_{N}|}\sum_{n\in\Phi_{N}}
f\otimes\bar{f}(\pi^{\times 2}((ca)^{n},(cb)^{n}))$, $a,b\in L$.

For $a,b\in L$,
the sequence $u_{n}=\pi^{\times2}((ca)^{n},(cb)^{n})$, $n\in\Z$, 
is contained in $X\times_{\hX}X=\bigl\{(x,y):\del(x)=\del(y)\bigr\}$,
where $\del$ is the factor mapping $X\ra\hX$.
If this sequence is well distributed in $X\times_{\hX}X$,
then $F(a,b)=\int_{X\times_{\hX}X}f\otimes\bar{f}\,d\mu_{X\times_{\hX}X}
=\int_{\hX}E(f|\hX)E(\bar{f}|\hX)\,d\mu_{\hX}=0$.
So, $F(a,b)\neq 0$
only if the sequence $(u_{n})$ is not well distributed in $X\times_{\hX}X$,
which only happens if the image $(\tu_{n})$ of $(u_{n})$ is not well distributed 
in the nil-maximal factor-torus $T$ of $X\times_{\hX}X$.
Using additive notation on $T$ we have $\tu_{n}=n\tc+n\ta+n\tb$, $n\in\Z$,
where $T$ contains the direct sum $S\oplus S$ of two copies of a torus $S$,
$\ta\in S\oplus\{0\}$, $\tb\in\{0\}\oplus S$,
and the sequence $(n\tc)$ is dense in the factor-torus $T/(S\oplus S)$.
The sequence $(\tu_{n})$ is not dense in $T$
only if the point $(\ta,\tb)$ is contained in a proper subtorus $R$ of $S\oplus S$,
and either for each $\tx\in S$ the fiber $\{\ty\in S:(\tx,\ty)\in R\}$ is a proper subtorus of $S$,
or for each $\ty\in S$ the fiber $\{\tx\in S:(\tx,\ty)\in R\}$ is a proper subtorus of $S$
(or both).
Without loss of generality, assume that the first possibility holds.
Then, returning back to $X\times_{\hX}X$,
we obtain that the sequence $(u_{n})$ is not well distributed in this space
only if the point $\pi^{\times2}(a,b)$ is contained in a subnilmanifold $D$
(the preimage of the torus $R$) in $Z\times Z$
with $D\ni 1_{X\times_{\hX}X}$
such that for every $x\in Z$ the fiber $W_{x}=\{y\in Z:(x,y)\in D\}$
is a proper normal subnilmanifold of $Z$.
Since, by our assumption, $\tro(W_{x})=0$ for all $x$,
we have $\tro^{\times2}(D)=0$,
so $\ro^{\times2}((\pi^{\times2})^{-1}(D))=0$.
The function $F(a,b)$ may only be nonzero
on the union of a countable collection of the subnilmanifolds $D$ like this,
so $\int_{L\times L}F(a,b)\,d\ro^{\times2}(a,b)=0$.
Hence, $\lim_{N\ras\infty}\frac{1}{|\Phi_{N}|}\sum_{n\in\Phi_{N}}|\phi(n)|^{2}=0$,
which means that $\phi\in\cZ$.
\endproof

Since, in the notation of \rfr{Lemma}{4.6}, 
the sequence $\hf(\pi(c^{n}))=\hf(c^{n}1_{X})$, $n\in\Z$, is a basic nilsequence,
$\phi$ is a sum of a nilsequence and a null-sequence,
so $\phi\in\cMb$ in this case as well.%
\endproof

The proof of \rfr{Proposition}{4.2} in the case $d\geq 2$ 
is not much harder than in the case $d=1$,
and we will only sketch it.
For each $\om\in\Om$, instead of a single element $a_{\om}\in G_{\om}$ 
we now have $d$ commuting elements $a_{\om,1},\ld,a_{\om,d}\in G_{\om}$.
After passing to a single nilmanifold $X=G/\Gam$,
we obtain $d$ mappings $\tau_{i}\col\Om\ra G$, $\om\mapsto a_{\om,i}$, $i=1,\ld,d$,
and so, the mapping $\tau=(\tau_{1},\ld,\tau_{d})\col\Om\ra G^{d}$.
We define a measure $\ro$ on $G^{d}$ by $\ro=\tau_{*}(\theta\nu)$;
then \rfr{Proposition}{4.2} follows from the following modification of \rfr{Proposition}{4.3}:
\proposition{P-mullimMes}{4.7}{}
Let $X=G/\Gam$ be a nilmanifold,
let $\ro$ be a finite Borel measure on $G^{d}$,
and let $f\in C(X)$.
Then the sequence  
$\phi(n_{1},\ld,n_{d})=\int_{G^{d}}f(\pi(a_{1}^{n_{1}}\ld a_{d}^{n_{d}}))\,d\ro(a_{1},\ld,a_{d})$,
$(n_{1},\ld,n_{d})\in\Z^{d}$, is in $\cMb$.
(If $X$ is an $r$-step nilmanifold, then $\phi\in\cMb_{r}$.)
\endproposition
\npar
The proof of this proposition is the same as of \rfr{Proposition}{4.3},
with $a^{n}$ replaced by $a_{1}^{n_{1}}\ld a_{d}^{n_{d}}$,
and the mapping $G\ra X$, $a\mapsto\pi(a)$, 
replaced by the mapping $G^{d}\ra X$, $(a_{1},\ld,a_{d})\mapsto\pi(a_{1}\ld a_{d})$.

\vbreak{-9000}{2mm}
Uniting \rfr{Proposition}{4.2} with \rfr{Lemma}{4.1},
we obtain: 
\theorem{P-IntM}{4.8}{}
Let $(\Om,\nu)$ be a measure space
and let $\Om\ra\cM$, $\om\mapsto\phi_{\om}$, be an integrable mapping.
Then the sequence $\phi(n)=\int_{\Om}\phi_{\om}(n)\,d\nu$ is in $\cM$ as well.
If, for some $r$, $\phi_{\om}\in\cM_{r}$ for all $\om$, then $\phi\in\cM_{r}$.
\endtheorem
\section{5}{Multiple polynomial correlation sequences and nilsequences}

Now let $(W,\B,\mu)$ be a probability measure space
and let $T$ be an ergodic invertible measure preserving transformation of $W$.
Let $p_{1},\ld,p_{k}$ be polynomials $\Z^{d}\ra\Z$.
Let $A_{1},\ld,A_{k}\in\B$ and let 
$\phi(n)=\mu\bigl(T^{p_{1}(n)}A_{1}\cap\ld\cap T^{p_{k}(n)}A_{k}\bigr)$,
$n\in\Z^{d}$;
or, more generally, let $f_{1},\ld,f_{k}\in L^{\infty}(W)$ and 
$\phi(n)=\int_{W}T^{p_{1}(n)}f_{1}\cd\ld\cd T^{p_{k}(n)}f_{k}d\mu$, 
$n\in\Z^{d}$.
Then, given $\eps>0$,
there exist an $r$-step nilsystem $(X,a)$, $X=G/\Gam$, $a\in G$,
and functions $\tf_{1},\ld,\tf_{k}\in L^{\infty}(X)$ such that, 
for $\ophi(n)=\int_{X}a^{p_{1}(n)}\tf_{1}\cd\ld\cd a^{p_{k}(n)}\tf_{k}d\mu_{X}$, $n\in\Z^{d}$,
the set $\bigl\{n\in\Z^{d}:|\ophi(n)-\phi(n)|>\eps\bigr\}$ has zero uniform density.
Moreover, there is a universal integer $r$
that works for all systems $(W,\B,\mu,T)$, functions $h_{i}$, and $\eps$,
and depends only on the polynomials $p_{i}$;
for the minimal such $r$,
the integer $c=r-1$ is called {\it the complexity\/} 
of the system $\{p_{1},\ld,p_{k}\}$ (see \brfr{L5}).

(Here is a sketch of the proof, for completeness;
for more details see \brfr{HK1} and \brfr{BHK}.
By \brfr{L3}, there exists $c\in\N$, which only depends on the polynomials $p_{i}$,
such that, if $(V,\nu,S)$ is an ergodic probability measure preserving system
and $Z_{c}(V)$ is the $c$-th Host-Kra-Ziegler factor of $V$
and $h_{1},\ld,h_{k}\in L^{\infty}(V)$ are such that $E(h_{i}|Z_{c}(V))=0$ for some $i$,
then for any F{\o}lner sequence $(\Phi_{N})$ in $\Z^{d}$ one has
$\lim_{N\ras\infty}\frac{1}{|\Phi_{N}|}\sum_{n\in\Phi_{N}}
\int_{V}S^{p_{1}(n)}h_{1}\cd\ld\cd S^{p_{k}(n)}h_{k}d\nu=0$.
Applying this to the ergodic components of the system $(W\times W,\mu\times\mu,T\times T)$
and the functions $h_{i}=f_{i}\tens\bar{f}_{i}$, $i=1,\ld,k$,
we obtain that for any F{\o}lner sequence $(\Phi_{N})$ in $\Z^{d}$,
\lequ{
\lim_{N\ras\infty}\frac{1}{|\Phi_{N}|}\sum_{n\in\Phi_{N}}
\Bigl|\int_{W\times W}T^{p_{1}(n)}f_{1}\cd\ld\cd T^{p_{k}(n)}f_{k}d\mu\Bigr|^{2}
\-\\\-
=\lim_{N\ras\infty}\frac{1}{|\Phi_{N}|}\sum_{n\in\Phi_{N}}
\int_{W\times W}T^{p_{1}(n)}f_{1}(x)\cd T^{p_{1}(n)}\bar{f}_{1}(y)\cd\ld
\-\\\-
\cd T^{p_{k}(n)}f_{k}(x)\cd T^{p_{k}(n)}\bar{f}_{k}(y)d(\mu(x)\times\mu(y))=0
}
whenever, for some $i$, the function $f_{i}\tens\bar{f_{i}}$
has zero conditional expectation 
with respect to almost all ergodic components of $Z_{c}(W\times W)$.
This is so if $E(f_{i}|Z_{c+1}(W))=0$,
and we obtain that the sequence $\int_{W}T^{p_{1}(n)}f_{1}\cd\ld\cd T^{p_{k}(n)}f_{k}d\mu$
tends to zero in uniform density whenever $E(f_{i}|Z_{r}(W))=0$ for some $i$, where $r=c+1$.
It follows that for any $f_{1},\ld,f_{k}\in L^{\infty}(W)$ the sequence
$$
\int_{W}T^{p_{1}(n)}f_{1}\cd\ld\cd T^{p_{k}(n)}f_{k}d\mu
-\int_{Z_{r}(W)}T^{p_{1}(n)}E(f_{1}|Z_{r}(W))\cd\ld\cd T^{p_{k}(n)}E(f_{k}|Z_{r}(W))d\mu_{Z_{r}(W)}
$$
tends to zero in uniform density.
Now, $Z_{r}(W)$ has the structure of the inverse limit of a sequence of $r$-step nilmanifolds
on which $T$ acts as a translation;
given $\eps>0$, we can therefore find an $r$-step nilmanifold factor $X$ of $W$
such that $\|E(f_{i}|Z_{r}(W))-E(f_{i}|X)\|_{L^{\infty}(W)}
<\eps/\prod_{j=1}^{k}\|f_{j}\|_{L^{\infty}(W)}$ for all $i$.
Putting $\tf_{i}=E(f_{i}|X)$, $i=1,\ld,k$,
and denoting the translation induced by $T$ on $X$ by $a$,
we then have
\lequ{
\Bigl|\int_{Z_{r}(W)}T^{p_{1}(n)}E(f_{1}|Z_{r}(W))\cd\ld\cd T^{p_{k}(n)}E(f_{k}|Z_{r}(W))d\mu_{Z_{r}(W)}
\-\\\-
-\int_{X}a^{p_{1}(n)}\tf_{1}\cd\ld\cd a^{p_{k}(n)}\tf_{k}d\mu_{X}\Bigr|
<\eps
} 
for all $n$,
which implies the assertion.)

So, there exists $\lam\in\cZ$ such that $\|\phi-(\psi+\lam)\|<\eps$.
After replacing $\tf_{i}$ by $L^{1}$-close continuous functions,
we may assume that $\tf_{1},\ld,\tf_{k}\in C(X)$,
and still $\|\phi-(\psi+\lam)\|<\eps$.
Applying \rfr{Theorem}{3.3} to the nilmanifold $X^{k}=G^{k}/\Gam^{k}$, 
the diagonal subnilmanifold $Y=\bigl\{(x,\ld,x),\ x\in X\bigr\}\sle X^{k}$,
the polynomial sequence $g(n)=(a^{p_{1}(n)},\ld,a^{p_{k}(n)})$, $n\in\Z^{d}$, in $G^{k}$, 
and the function $f(x_{1},\ld,x_{k})=\tf_{1}(x_{1})\cd\ld\cd\tf_{k}(x_{k})\in C(X^{k})$,
we obtain that $\psi\in\cMb(r,s)$, so also $\psi+\lam\in\cMb(r,s)$.
Since $\eps$ is arbitrary and, by \rfr{Lemma}{2.3}, $\cM_{r,s}$ is the closure of $\cMb_{r,s}$,
we obtain:
\proposition{P-ergCor}{5.1}{}
Let $(W,\B,\mu,T)$ be an ergodic invertible probability measure preserving system,
let $f_{1},\ld,f_{k}\in L^{\infty}(W)$,
and let $p_{1},\ld,p_{k}$ be polynomials $\Z^{d}\ra\Z$.
Then the sequence 
$\phi(n)=\int_{W}T^{p_{1}(n)}f_{1}\cd\ld\cd T^{p_{k}(n)}f_{k}d\mu$, $n\in\Z^{d}$, 
is in $\cM$.
If the complexity of the system $\{p_{1},\ld,p_{k}\}$ is $c$ and $\deg p_{i}\leq s$ for all $i$,
then $\phi_{n}\in\cM_{c+1,s}$.
\endproposition

Let now $(W,\B,\mu,T)$ be a non-ergodic (or, rather, not necessarily ergodic) system.
Let $\mu=\int_{\Om}\mu_{\om}d\nu(\om)$ be the ergodic decomposition of $\mu$.
For each $\om\in\Om$, 
let $\phi_{\om}(n)=\int_{W}T^{p_{1}(n)}f_{1}\cd\ld\cd T^{p_{k}(n)}f_{k}d\mu_{\om}$, $n\in\Z^{d}$;
then $\om\mapsto\phi_{\om}$ is a measurable mapping $\Om\ra\Li$,
and $\phi(n)=\int_{\Om}\phi_{\om}(n)\,d\nu(\om)$, $n\in\Z^{d}$.
By \rfr{Proposition}{5.1}, 
for each $\om\in\Om$ we have $\phi_{\om}\in\cM_{c+1,s}\sle\cM_{l}$, where $l=2(c+1)s$.
By \rfr{Theorem}{4.8} we obtain:
\theorem{P-Cor}{5.2}{}
Let $(W,\B,\mu,T)$ be an invertible probability measure preserving system,
let $f_{1},\ld,f_{k}\in L^{\infty}(W)$,
and let $p_{1},\ld,p_{k}$ be polynomials $\Z^{d}\ra\Z$.
Then the sequence
$\phi(n)=\int_{W}T^{p_{1}(n)}f_{1}\cd\ld\cd T^{p_{k}(n)}f_{k}d\mu$, $n\in\Z^{d}$,
is in $\cM$.
If the complexity of the system $\{p_{1},\ld,p_{k}\}$ is $c$ and $\deg p_{i}\leq s$ for all $i$,
then $\phi_{n}\in\cM_{l}$, where $l=2(c+1)s$.
\endtheorem

Since $\cM\sle\cP+\cZ$, 
where $\cP$ is the closure in $\Li$ of the algebra of bounded generalized polynomials
(see the last paragraph of \rfr{Section}{2}),
we get as a corollary:
\corollary{P-GPS}{5.3}{}
Up to a null-sequence, 
the sequence $\phi$ is uniformly approximable by generalized polynomials.
\endcorollary
\bibliography{}
\bibart BHK
a:V. Bergelson, B. Host and B. Kra
t:Multiple recurrence and nilsequences.
j:Inventiones Math.
n:160
y:2005
p:\no{2}, 261--303
*
\bibart BL
a:V. Bergelson and A. Leibman
t:Distribution of values of bounded generalized polynomials
j:Acta Math.
n:198
y:2007
p:155-230
*
\bibook F
a:H. Furstenberg
t:Recurrence in Ergodic Theory and Combinatorial Number Theory
i:Princeton Univ. Press, 1981
*
\bibart GT
a:B. Green and T. Tao
t:The quantitative behaviour of polynomial orbits on nilmanifolds
j:Annals of Math.
n:175
y:2012
p:\no{2}, 465-540
*
\bibart H{\aa}1
a:I.J. H{\aa}land (Knutson)
t:Uniform distribution of generalized polynomials
j:J. Number Theory
n:45
y:1993 
p:327-366 
*
\bibart H{\aa}2
a:I.J. H{\aa}land (Knutson)
t:Uniform distribution of generalized polynomials of the product type
j:Acta Arith. 
n:67
y:1994 
p:13-27 
*
\bibart HK1
a:B. Host and B. Kra 
t:Non-conventional ergodic averages and nilmanifolds
j:Annals of Math.
n:161
y:2005
p:\no{1}, 397-488
* 
\bibart HK2
a:B. Host and B. Kra 
t:Nil-Bohr sets of integers
j:Ergodic Theory and Dynam. Systems
n:31
y:2011
p:113-142
* 
\bibart HKM
a:B. Host, B.Kra, and A. Maass
t:Nilsequences and a structure theorem for topological dynamical systems
j:Adv. in Math. 
n:224 
y:2010 
p:103-129
*
\bibart L1
a:A. Leibman
t:Pointwise convergence of ergodic averages 
for polynomial sequences of translations on a nilmanifold
j:Ergodic Theory and Dynam. Systems
n:25
y:2005
p:201-213
*
\bibart L2
a:A. Leibman
t:Pointwise convergence of ergodic averages 
for polynomial actions of $\Z^{d}$ by translations on a nilmanifold
j:Ergodic Theory and Dynam. Systems
n:25
y:2005
p:215-225
*
\bibart L3
a:A. Leibman
t:Convergence of multiple ergodic averages along polynomials of several variables
j:Israel J. of Math.
n:146
y:2005
p:303-315
*
\bibart L4
a:A. Leibman
t:Rational sub-nilmanifolds of a compact nilmanifold
j:Ergodic Theory and Dynam. Systems
n:26
y:2006
p:787-798
*
\bibart L5
a:A. Leibman
t:Orbit of the diagonal in the power of a nilmanifold
j:Trans. AMS
n:362
y:2010
p:1619-1658
*
\bibart L6
a:A. Leibman
t:Multiple polynomial correlation sequences and nilsequences
j:Ergodic Theory and Dynam. Systems
n:30
y:2010
p:841-854
*
\bibart L7
a:A. Leibman
t:A canonical form and the distribution of values of generalized polynomials 
j:Israel J. Math
n:188
y:2012
p:131-176
*
\bibart M
a:A. Malcev
t:On a class of homogeneous spaces
j:Amer. Math. Soc. Transl.
n:9
y:1962
p:276-307
*
\bibart Z
a:T. Ziegler
t:Universal characteristic factors and Furstenberg averages
j:J. Amer. Math. Soc.
n:20
y:2007
p:\no{1}, 53-97
* 
\endbibliography
\bye